\documentclass[10pt, a4paper]{article}
\usepackage{graphics}
\usepackage{latexsym}
\usepackage{amssymb}
\usepackage{amsmath}
\usepackage{color}
\usepackage{epsfig}

\newcommand{\beq}{\begin{equation}}
\newcommand{\eeq}{\end{equation}}
\newcommand{\beqas}{\begin{eqnarray*}}
\newcommand{\eeqas}{\end{eqnarray*}}
\newcommand{\ep}{\varepsilon}

\newcommand{\ue}{u^{\varepsilon}}

\newcommand{\bproof}{{\it Proof}\ \ }
\newcommand{\eproof}{{\hfill$\Box$}}

\newcommand{\be}{\begin{equation}}
\newcommand{\ee}{\end{equation}}
\newcommand{\bea}{\begin{eqnarray}}
\newcommand{\eea}{\end{eqnarray}}
\newcommand{\beas}{\begin{eqnarray*}}
\newcommand{\eeas}{\end{eqnarray*}}

\def\IR{\mathbb{R}}

\def\curl{{\rm curl\,}}

\def\scurl{{\rm curl}}

\newtheorem{theorem}{Theorem}[section]
\newtheorem{lemma}[theorem]{Lemma}

\newtheorem{proposition}[theorem]{Proposition}
\newtheorem{definition}[theorem]{Definition}
\newtheorem{remark}[theorem]{Remark}

\definecolor{darkgreen}{rgb}{0,0.7,0}

\numberwithin{equation}{section}
\numberwithin{equation}{section}
\title{Homogenization error for two scale Maxwell equations
\thanks{%
VTC is supported   by an NTU graduate scholarship, VHH is supported by the  AcRF Tier 1 research grant RG69/10, the Singapore A*Star SERC grant 122-PSF-0007 and the AcRF Tier 2 grant MOE 2013-T2-1-095, ARC 44/13. 
}}

\author{
Van Tiep Chu
and 
Viet Ha Hoang
\thanks{Division of Mathematical Sciences, 
      School of Physical and Mathematical Sciences, 
       Nanyang Technological University, Singapore 637371 
       }
}
\date{}
\begin{document}

\maketitle

\begin{abstract}
For two scale elliptic equations in a domain $D$, standard homogenization errors are deduced with the assumption that the solution $u_0$ of the homogenized equation belongs to $H^2(D)$. For two scale Maxwell equations, the corresponding required regularity is $u_0\in H^1(\curl, D)$. These regularity conditions normally do not hold in general polygonal domains, which are of interests for finite element discretization. The paper establishes homogenization errors when $u_0$ belongs to a weaker regularity space $H^{1+s}(D)$ for elliptic problems and $H^s(\curl,D)$ for Maxwell problems where $0<s<1$. Though we only present the results for two scale Maxwell equations when $u_0\in H^s(\curl,D)$ with $0<s<1$, the procedure works verbatim for elliptic equations when $u_0$ belongs to $H^{1+s}(D)$ with $0<s<1$. 
\end{abstract}

\section{Introduction}
For two scale elliptic problems in a domain $D\subset\IR^d$, standard homogenization error estimates are deduced with the assumption that the solution $u_0$ of the homogenized equation belongs to $H^2(D)$ (see, e.g., \cite{BLP} and \cite{JKO}). In many cases, this condition does not hold. A typical situation is that of nonconvex polygonal domains which we often need to consider in the context of finite element discretization. For Maxwell equations, the approaches in \cite{BLP}, \cite{JKO} and \cite{CD}  require the regularity $u_0\in H^1(\curl, D)$ (see Section 3) which normally does not hold in polygonal domains. However, in polygonal domains, the solutions $u_0$ belong to a weaker regularity space, i.e $u_0\in H^{1+s}(D)$ for elliptic problems and $u_0\in H^s(\curl,D)$ for Maxwell problems with $0<s<1$ (see, e.g, \cite{Hiptmair}). The paper develops new homogenization errors for two scale Maxwell problems where $u_0\in H^s(\curl,D)$ for $0<s<1$. Though we only present the results for two scale Maxwell equations, our approach applies verbatim for two scale elliptic equations with $u_0\in H^{1+s}(D)$ for $0<s<1$. 

The paper is organized as follows. In the next section, we introduce the two scale Maxwell equation \eqref{eq:msprob} and consider its homogenization. Homogenization limit of Maxwell equation \eqref{eq:msprob} is deduced in Benssousan et al. \cite{BLP} though the two scale asymptotic expansion is performed only for the case where the coefficient $b^\ep$ in \eqref{eq:msprob} is a constant isotropic matrix; and as a consequence the corrector function ${\frak u}_1$ in \eqref{eq:fu1} is not derived. To get the correctors explicitly, we therefore employ the two scale convergence method (\cite{Nguetseng} and \cite{Allaire}) to establish the two scale homogenized equation, from which we deduce explicit formulae for the two correctors $u_1$ and ${\frak u}_1$ in \eqref{eq:u1} and \eqref{eq:fu1} respectively. We note that two scale convergence limits of bounded sequences in $H(\curl,D)$ are considered in \cite{WK} but the form of the limit function is slightly different from the form that we need in this paper. We therefore derive and prove in details the two scale convergence limits. Section 3 establishes the homogenization error when $u_0\in H^1(\curl,D)$. We follow the standard approach of \cite{JKO} for elliptic equations, but it appears that this has not been done for Maxwell equations in the literature. Our main contribution is contained in Section 4 where we derive homogenization errors for the case where $u_0\in H^s(\curl,D)$ for $0<s<1$. The last section proves the regularity required for the solutions of the cell problems \eqref{eq:cellprobw} and \eqref{eq:cellprobN}, and for $u_0$. 

Throughout the paper, by $\nabla$ and $\curl$, without indicating explicitly the variable, we denote the total gradient and $\curl$ of a function of the variable $x$. Partial gradient and partial $\curl$ of a function that depends on $x$ and $y$ are denoted by $\nabla_x$, $\nabla_y$, $\scurl_x$ and $\scurl_y$ respectively. Repeated indices indicate summation. The notation $\#$ denotes spaces of periodic functionsd; $c$ denotes various constants of different values. 

\section{Homogenization of the Maxwell equation}
Let $D$ be a domain in $\IR^3$ with a Lipschitz boundary. Let $Y$ be the open unit cube in $\IR^3$. Let $a(x,y), b(x,y): D\times Y\to\IR^{3\times 3}$ be positive definite matrix functions, i.e. there are constants $c_1, c_2>0$ so that for all vectors $\xi,\zeta\in\IR^3$ we have
\be
c_1|\xi|^2\le a(x,y)\xi\cdot\xi,\ \ \ c_1|\xi|^2\le b(x,y)\xi\cdot\xi,
\label{eq:coercive}
\ee
and 
\be
a(x,y)\xi\cdot\zeta\le c_2|\xi||\zeta|,\ \ \ b(x,y)\xi\cdot\zeta\le c_2|\xi||\zeta|
\label{eq:bounded}
\ee
for all $x\in D$ and $y\in Y$, where $|\cdot|$ denotes the Euclid norm in $\IR^3$. We define the multiscale coefficients $a^\ep:D\to\IR^{3\times 3}$ and $b^\ep:D\to\IR^{3\times 3}$ as 
\[
a^\ep(x)=a\left(x,{x\over\ep}\right),\ \ \ b^\ep(x)=b\left(x,{x\over\ep}\right).
\]
Let $V=H_0(\curl,D)$ and $f\in V'$. 
We consider the problem
\be
\curl(a^\ep(x)\curl\ue)+b^\ep\ue=f,\ \ \ \ue\times n=0,\ \ \mbox{on}\ \partial D,
\label{eq:msprob}
\ee
where $n$ is the outward normal vector on $\partial D$. In variational form, this problem is written as
\be
\int_D[a^\ep(x)\curl\ue\cdot\curl\phi+b^\ep\ue\cdot\phi]dx=\int_Df\cdot\phi dx\ \ \forall\phi\in V.
\label{eq:msprobweak}
\ee

We study homogenization of equation \eqref{eq:msprob} via two scale convergence in this section. We first recall the concept of two scale convergence (see \cite{Nguetseng} and \cite{Allaire}).
\begin{definition}
A sequence of functions $\{w^\ep\}_\ep\subset L^2(D)$ two scale converges to a function $w_0\in L^2(D\times Y)$ if for all smooth functions $\phi\in C^\infty(D\times Y)$ which are periodic with respect to $y$ with the period being $Y$
\[
\lim_{\ep\to 0}\int_Dw^\ep(x)\phi(x,{x\over\ep})dx=\int_D\int_Y w_0(x,y)\phi(x,y)dy dx.
\]
\end{definition}
The definition makes sense due to the following result.
\begin{proposition}
From a bounded sequence in $L^2(D)$ we can extract a two scale convergent subsequence.
\end{proposition}
To deduce the homogenized problem for the two scale Maxwell equation \eqref{eq:msprob} by two scale convergence, we first establish two scale convergence results for a bounded sequence in $H(\curl,D)$. These results are first established in Wellander and Kristensson \cite{WK}. However, for completeness we present the proof here as the convergence limit we use is in a slightly different form from that of \cite{WK} (see Remark \ref{rem:WK}).

Let $H_\#(\curl,Y)$ be the space of functions in $H(\curl,Y)$ which can be extended periodically to $H_{\rm loc}(\curl,\IR^3)$. This space is the closure of the space $C_\#^\infty(Y)$ of periodic smooth functions in $H(\curl,Y)$. 
We denote by  $\tilde H_\#(\curl,Y)$ the space of equivalent classes of functions in $H_\#(\curl,Y)$ where $\phi$ and $\psi$ are regarded as equal if $\curl_y\phi=\curl_y\psi$, with the norm $\|\phi\|_{\tilde H_\#(\curl,Y)}=\|\scurl_y\phi\|_{L^2(Y)^3}$.  
\begin{proposition}
Let $\{w^\ep\}_\ep\subset H(\curl,D)$. There is a subsequence (not renumbered), and functions $w_0\in L^2(D)^3$, $\frak w_1\in L^2(D,H^1_\#(Y)/\IR)$, $w_1\in L^2(D,\tilde H_\#(\curl,Y))$ such that
\[
w^\ep \stackrel{\mbox{two scale}}{\longrightarrow}w_0+\nabla_y{\frak w}_1,\ \ \mbox{and}\ \ \curl w^\ep \stackrel{\mbox{two scale}}{\longrightarrow}\curl w_0+\scurl_y w_1.
\]
\end{proposition}
\bproof
As $\{w^\ep\}_\ep$ is bounded in $H(\curl,D)$, there is a subsequence (not renumbered) such that $w^\ep$ and $\curl w^\ep$ are two scale convergent. 
Let $\xi\in L^2(D\times Y)^3$ be the two scale limit of $\{w^\ep\}_\ep$. Let $\phi(x,y)=\ep\Phi(x,y)$ where $\Phi$ is a smooth function in $C_0^\infty(D,C^\infty_\#(Y))^3$. We have
\beqas
0=\lim_{\ep\to 0}\int_D\curl w^\ep\cdot\ep\Phi(x,{x\over\ep})dx=\lim_{\ep\to 0}\int_D w^\ep\cdot\ep\curl\Phi(x,{x\over\ep})dx\\
=\lim_{\ep\to 0}\int_D w^\ep\cdot\scurl_y\Phi(x,{x\over\ep})dx\\
=\int_D\int_Y\xi(x,y)\cdot\scurl_y\Phi(x,y)dydx.
\eeqas
As this holds for all the smooth functions $\Phi(x,y)\in C_0^\infty(D,C_\#^\infty(Y))$, we have that
\[
\xi(x,y)=\xi_0(x)+\nabla_y{\frak w}_1(x,y)
\]
for a function ${\frak w}_1(x,y)\in L^2(D, H^1_\#(Y)/\IR)$ and a function $\xi_0\in L^2(D)^3$. As $\int_Y\xi(x,y)dy=\xi_0(x)$, $\xi_0$ equals the weak limit $w_0$ of $w^\ep$ in $L^2(D)^3$. 

Let $\eta(x,y)\in L^2(D\times Y)^3$ be the two scale limit of $\curl w^\ep$. For $\Phi\in C^\infty_0(D,C^\infty_\#(Y))$, we have
\[
\int_D\curl w^\ep\cdot\nabla\Phi(x,{x\over\ep})dx=\int_D w^\ep\cdot\curl\nabla\Phi(x,{x\over\ep})dx-\int_{\partial D}(w^\ep\times n)\cdot\nabla\Phi(x,{x\over\ep})ds=0.
\]
Thus
\beqas
0=\lim_{\ep\to 0}\int_D\curl w^\ep\cdot\ep\nabla\Phi(x,{x\over\ep})=\lim_{\ep\to 0}\int_D\curl w^\ep\cdot\nabla_y\Phi(x,{x\over\ep})dx\\
=\int_D\int_Y\eta(x,y)\cdot\nabla_y\Phi(x,y)dydx.
\eeqas
Thus there is a function $w_1(x,y)\in L^2(D,\tilde H_\#(\curl, Y))$ and a function $\eta_0(x)\in L^2(D)^3$ such that 
\[
\eta(x,y)=\eta_0(x)+\scurl_yw_1(x,y).
\]
As $\int_Y\eta(x,y)dy=\eta_0(x)$, $\eta_0$ equals the weak limit of $\curl w^\ep$ in $L^2(D)^3$, which is $\curl w_0$. 
\eproof
\begin{remark}\label{rem:WK} Wellander and Kristensson \cite{WK} show that if $w^\ep$ two scale converges to $w(x,y)$, then there is a function $\tilde w_1(x,y)\in L^2(D,H_\#(\curl, D))$ so that the two scale limit $\eta(x,y)$ of $\curl w^\ep$ satisfies $\eta(x,y)=\scurl_xw(x,y)+\scurl_y\tilde w_1(x,y)$. As $w(x,y)=w_0(x)+\nabla_y \frak w_1(x,y)$, $\scurl_xw(x,y)=\curl w_0(x)+\scurl_x\nabla_y\frak w_1(x,y)=\curl w_0(x)-\scurl_y\nabla_x\frak w_1(x,y)$. Therefore 
\[
\eta(x,y)=\curl w_0(x)+\scurl_y(-\nabla_x\frak w_1(x,y)+\tilde w_1(x,y)).
\]
Our result above is consistent with the result by Wellander and Kristensson \cite{WK}. 
\end{remark}

We now derive the two scale homogenized problem for the two scale Maxwell equation \eqref{eq:msprob}. 

\begin{proposition}
There are functions $u_0(x)\in H_0(\curl, D)$, $u_1(x,y)\in L^2(D,\tilde H_\#(\curl,Y))$ and ${\frak u_1}(x,y)\in L^2(D, H^1_\#(Y)/\IR)$ such that
\[
\ue\stackrel{\mbox{two scale}}{\longrightarrow}u_0+\nabla_y{\frak u_1},
\]
and 
\[
\curl u^\ep\stackrel{\mbox{two scale}}{\longrightarrow}\curl u_0+\curl_yu_1.
\]
The functions $u_0$, $u_1$ and ${\frak u}_1$ satisfy the problem
\begin{eqnarray}
\int_D\int_Y\left[a(x,y)(\curl u_0+\scurl_yu_1)\cdot(\curl v_0+\scurl_yv_1)+b(x,y)(u_0+\nabla_y\frak u_1)\cdot(v_0+\nabla_y\frak v_1)\right]dydx\nonumber\\
=\int_Df(x)\cdot v_0(x)dx
\label{eq:2shomprob}
\end{eqnarray}
for all $v_0\in H_0(\curl, D)$, $v_1\in L^2(D,\tilde H_\#(\curl, Y))$ and ${\frak v}_1\in L^2(D,H^1_\#(Y)/\IR)$. 
\end{proposition}
\bproof
The proof is standard. From \eqref{eq:coercive} and \eqref{eq:bounded}, $\{\ue\}_\ep$ is uniformly bounded in $H(\curl,D)$. Thus we can extract a subsequence (not renumbered) so that 
\[
\ue\stackrel{\mbox{two scale}}{\longrightarrow}u_0+\nabla_y{\frak u_1},
\ \ \ 
\mbox{and}\ \  
\curl u^\ep\stackrel{\mbox{two scale}}{\longrightarrow}\curl u_0+\curl_yu_1
\]
for $u_0(x)\in H_0(\curl, D)$, $u_1(x,y)\in L^2(D,\tilde H_\#(\curl,Y))$ and ${\frak u_1}(x,y)\in L^2(D, H^1_\#(Y)/\IR)$.
We choose  functions $v_0\in C_0^\infty(D)^3$, $v_1\in C_0^\infty(D,C^\infty_\#(Y))^3$ and ${\frak v}_1\in C^\infty_0(D,C^\infty_\#(Y))$. Let 
\[
\phi(x)=v_0(x)+\ep v_1(x,{x\over\ep})+\ep\nabla{\frak v}_1(x,{x\over\ep}).
\]
as the test function in \eqref{eq:msprobweak}. We get 
\beqas
\int_D\left[a(x,{x\over\ep})\curl u^\ep(x)\cdot\left(\curl v_0(x)+\ep\scurl_xv_1(x,{x\over\ep})+\scurl_yv_1(x,{x\over\ep})\right)\right.\\
+\left.b(x,{x\over\ep})u^\ep(x)\cdot \left(v_0(x)+\ep v_1(x,{x\over\ep})+\ep\nabla_x{\frak v}_1(x,{x\over\ep})+\nabla_y{\frak v}(x,{x\over\ep})\right)\right]\\
=\int_Df(x)\cdot\left[ v_0(x)+\ep v_1(x,{x\over\ep})+\ep\nabla_x{\frak v}_1(x,{x\over\ep})+\nabla_y{\frak v}_1(x,{x\over\ep})\right].
\eeqas
Passing to the two scale limit, we have
\begin{eqnarray}
\int_D\int_Y\left[ a(x,y)(\curl u_0+\scurl_yu_1)\cdot(\curl v_0+\scurl_yv_1)+b(x,y)(u_0+\nabla_y{\frak u}_1)\cdot(v_0+\nabla_y{\frak v}_1)\right]dydx\nonumber\\
=\int_D f(x)\cdot v_0(x)dx+\int_D\int_Y f(x)\cdot\nabla_y{\frak v}_1(x,y)dydx=\int_D f(x)\cdot v_0(x)dx.
\label{eq:ms1}
\end{eqnarray}
Equation \eqref{eq:2shomprob} follows from a density argument. 
From Lax-Milgram Lemma, \eqref{eq:ms1} has a unique solution $(u_0,u_1,{\frak u}_1)$. The two scale convergence properties thus hold for the whole sequence $\{\ue\}_\ep$.
\eproof

We now derive the cell problems and the homogenized equation. 

Letting $v_0=0$, $v_1=0$, and deduce that 
\[
\int_D\int_Yb(x,y)(u_0+\nabla_y\frak u_1)\cdot \nabla_y\frak v_1dy dx=0.
\]
For each $r=1,2,3$, let $w^r(x,\cdot)\in L^2(D, H^1_\#(Y)/\IR)$ be the solution of the cell problem
\be
\int_D\int_Yb(x,y)(e_r+\nabla_yw^r)\cdot\nabla_y\psi dydx=0\ \ \forall\psi\in L^2(D,H^1_\#(Y)/\IR)
\label{eq:cellprobw}
\ee
where $e_r$ is the vector in $\IR^3$ with all the components being 0, except the $r$th component which equals 1. This is the standard cell problem in elliptic homogenization. From this we have 
\be
\frak u_1(x,y)=w^r(x,y)u_{0r}(x).
\label{eq:fu1}
\ee
 Therefore
\[
\int_D\int_Yb(x,y)(u_0+\nabla_y\frak u_1)\cdot v_0dxdy=\int_Db^0(x)u_0(x)\cdot v_0(x)dx,
\]
where the positive definite matrix $b^0(x)$ is defined as 
\be
b^0_{ij}(x)=\int_Y b(x,y)(e_j+\nabla w^j(x,y))\cdot(e_i+\nabla_yw^i(x,y))dy. 
\label{eq:b0}
\ee
Let $v_0=0$ and $\frak v_1=0$. We have 
\[
\int_D\int_Ya(x,y)(\curl u_0+\scurl_yu_1)\cdot\scurl_yv_1dydx=0
\]
for all $v_1\in L^2(D, \tilde H_\#(\curl,Y))$. For each $r=1,2,3$, let $N^r\in L^2(D,\tilde H_\#(\curl,Y))$ be the solution of
\be
\int_D\int_Ya(x,y)(e_r+\scurl_yN^r)\cdot\scurl_yvdydx=0
\label{eq:cellprobN}
\ee
for all $v\in L^2(D, \tilde H_\#(\curl,Y))$. We have
\be
u_1=\curl u_0(x)_rN^r(x,y).
\label{eq:u1}
\ee
The homogenized coefficient $a^0$ is determined by
\be
a^0_{ij}(x)=\int_Ya(x,y)(e_j+\scurl_yN^j)\cdot(e_i+\scurl_yN^i)dy.
\label{eq:a0}
\ee
We have
\[
\int_D\int_Ya(x,y)(\curl u_0+\scurl_yu_1)\cdot\curl v_0dxdy=\int_Da^0(x)\curl u_0(x)\cdot\curl v_0(x)dx\ \ \forall\,v_0\in H_0(\curl,D).
\]
The homogenized problem is
\be
\int_D [a^0(x)\curl u_0(x)\cdot\curl v_0(x)+b^0(x)u_0(x)\cdot v_0(x)]dx=\int_Df(x)\cdot v_0(x)dx\ \ \forall\,v_0\in H_0(\curl,D).
\label{eq:homprob}
\ee
\section{Homogenization error when $u_0\in H^1(\curl, D)$}
We show in this section the homogenization error when the solution $u_0$ of the homogenized problem \eqref{eq:homprob} belongs to the regularity space $H^1(\curl,D)$. We follow the approach in \cite{JKO}. 

\begin{theorem}\label{thm:main1}
Assume that $a(x,y)\in C(\bar D, C(\bar Y))$, $u_0\in H^1(\curl;D)$,  $N^r\in C^1(\bar D,C(\bar Y))^3$, $\scurl_yN^r\in C^1(\bar D,C(\bar Y))^3$, $w^r\in C^1(\bar D,C^1(\bar Y))$ \footnote{Indeed for Theorems \ref{thm:main1} and \ref{thm:main2}, we only need weaker regularity conditions $N^r\in W^{1,\infty}(D,L^\infty(Y))$ and $\scurl_yN^r\in W^{1,\infty}(D,L^\infty(Y))$, and $w^r\in W^{1,\infty}(D,W^{1,\infty}(Y))$.}for all $r=1,2,3$, then
\[
\|\ue-[u_0+\nabla_y{\frak u}_1(\cdot,{\cdot\over\ep})]\|_{L^2(D)^3}\le c\ep^{1/2}
\] 
and 
\[
\|\curl\ue-[\curl u_0+\curl_yu_1(\cdot,{\cdot\over\ep})]\|_{L^2(D)^3}\le c\ep^{1/2}.
\]
\end{theorem}
\bproof 
We consider the function 
\[
u_1^\ep(x)=u_0(x)+\ep N^r(x,{x\over\ep})\curl u_{0}(x)_r+\ep\nabla\frak u_1(x,{x\over\ep}).
\] 
We have
\beqas
\curl(a^\ep\curl u_1^\ep)+b^\ep u_1^\ep=\curl a(x,{x\over\ep})\Big[\curl u_0(x)+\ep\scurl_x N^r(x,{x\over\ep})\curl u_0(x)_r+\scurl_yN^r(x,{x\over\ep})\curl u_{0}(x)_r+\\
\ep \nabla \curl u_{0}(x)_r\times N^r(x,{x\over\ep})\Big]+\\
+b(x,{x\over\ep})\Big[u_0(x)+\nabla_yw^r(x,{x\over\ep})u_{0r}(x)+\ep N^r(x,{x\over\ep})\curl u_{0}(x)_r+\ep\nabla_x w^r(x,{x\over\ep})u_{0r}(x)+\ep w^r(x,{x\over\ep})\nabla u_{0r}(x)\Big]\\
=\curl a^0(x)\curl u_0(x)+b^0(x)u_0(x)+\curl G_r(x,{x\over\ep})\curl u_0(x)_r+g_r(x,{x\over\ep})u_{0r}(x)+\ep\curl I(x)+\ep J(x)
\eeqas
where the vector functions $G_r(x,y)$ and $g_r(x,y)$ are defined by
\be
(G_r)_i(x,y)=a_{ir}(x,y)+a_{ij}(x,y)\scurl_yN^r(x,y)_j-a^0_{ir}(x),
\label{eq:G}
\ee
\be
(g_r)_i(x,y)=b_{ir}(x,y)+b_{ij}(x,y){\partial w^r\over\partial y_j}(x,y)-b^0_{ir}(x);
\label{eq:g}
\ee
and
\beqas
I(x)=a(x,{x\over\ep})\Big[\curl_x N^r(x,{x\over\ep})\curl u_0(x)_r+\nabla\curl u_0(x)_r\times N^r(x,{x\over\ep})\Big]
\eeqas
and
\beqas
J(x)= 
b(x,{x\over\ep})\Big[N^r(x,{x\over\ep})\curl u_{0}(x)_r+\nabla_x w^r(x,{x\over\ep})u_{0r}(x)+w^r(x,{x\over\ep})\nabla u_{0r}(x)\Big].
\eeqas
From \eqref{eq:cellprobN}, we have that $\scurl_yG_r(x,y)=0$. Further from \eqref{eq:a0} $\int_YG_r(x,y)dy=0$. We thus deduce that there is a  function $\tilde G_r(x,y)$ such that $G_r(x,y)=\nabla_y \tilde G_r(x,y)$. From \eqref{eq:cellprobw}, we have that ${\rm div}_yg_r(y)=0$ and from \eqref{eq:b0} $\int_Yg_r(x,y)dy=0$. Therefore there is a function $\tilde g_r$ such that $g_r(x,y)=\scurl_y\tilde g_r(x,y)$. 
Thus for all $\phi\in{\cal D}(D)^3$
\beqas
&&\langle\curl(a^\ep\curl u_1^\ep)+b^\ep u_1^\ep-\curl (a^0\curl u_0)-b^0u_0,\phi\rangle=\\
&&\qquad\int_DG_r(x,{x\over\ep})\curl u_0(x)_r\cdot\curl\phi(x)dx+\int_D g_r(x,{x\over\ep})u_{0r}(x)\cdot\phi(x)dx\\
&&\qquad\qquad\qquad+\ep\int_DI(x)\cdot\curl\phi(x)dx+\ep\int_DJ(x)\cdot\phi(x)dx\\
&&\qquad=-\ep\int_D\tilde G_r(x,{x\over\ep}){\rm div}(\curl u_0(x)_r\curl\phi(x))dx-\ep\int_D\nabla_x\tilde G_r(x,{x\over\ep})\curl u_0(x)_r\cdot\curl\phi(x)dx\\
&&\qquad\qquad+\ep\int_D\tilde g_r(x,{x\over\ep})\cdot\curl(u_{0r}(x)\phi(x))dx-\ep\int_D\scurl_x\tilde g_r(x,{x\over\ep})u_{0r}(x)\cdot\phi(x)dx\\
&&\qquad\qquad+\ep\int_DI(x)\cdot\curl\phi(x)dx+\ep\int_DJ(x)\cdot\phi(x).
\eeqas
We note that
\[
\int_D\tilde G_r(x,{x\over\ep}){\rm div}(\curl u_0(x)_r\curl\phi(x))dx=\int_D\tilde G_r(x,{x\over\ep})\nabla\curl u_0(x)_r\cdot\curl\phi(x)dx
\]
and 
\[
\int_D\tilde g_r(x,{x\over\ep})\cdot\curl(u_{0r}(x)\phi(x))dx=\int_D\tilde g_r(x,{x\over\ep})\cdot(u_{0r}(x)\curl\phi(x)+\nabla u_{0r}(x)\times\phi(x))dx.
\]
As $\nabla_y\tilde G(x,\cdot)=G(x,\cdot)\in H^1(Y)^3$ so $\Delta_y\tilde G(x,\cdot)\in L^2(Y)$. Thus $\tilde G(x,\cdot)\in H^2(Y)$ which implies $\tilde G(x,\cdot)\in C(\bar Y)$. As $G(x,\cdot)\in C^1(\bar D,H^1_\#(Y)^3)$, we deduce that $\tilde G(x,y)\in C^1(\bar D, H^2(Y))\subset C^1(\bar D,C(\bar Y))$. The construction of $\tilde g_r$ in Jikov et al. \cite{JKO} implies that $g_r\in C^1(\bar D, C(\bar Y))$ (see  Hoang and Schwab \cite{HSmultirandom}). Thus
\[
\Big|\int_D\tilde G_r(x,{x\over\ep}){\rm div}(\curl u_0(x)_r\curl\phi(x))dx\le c\|\curl\phi\|_{L^2(D)^3},
\]
and 
\[
\Big|\int_D\tilde g_r(x,{x\over\ep})\cdot\curl(u_{0r}(x)\phi(x))dx\le c(\|\curl\phi\|_{L^2(D)^3}+\|\phi\|_{L^2(D)^3}).
\]
As $N^r\in C^1(\bar D,C(\bar Y))$ and $w^r\in C^1(\bar D,C(\bar Y))$, $\|I\|_{L^2(D)}$ and $\|J\|_{L^2(D)}$ are uniformly bounded with respect to $\ep$. 
From these we  conclude that 
\[
|\langle\curl(a^\ep\curl u_1^\ep)+b^\ep u_1^\ep-\curl(a^0\curl u_0)-b^0u_0,\phi\rangle|\le c\ep(\|\curl\phi\|_{L^2(D)^3}+\|\phi\|_{L^2(D)^3})
\]
where $\langle\cdot\rangle$ denotes the duality pairing of $V$ and $V'$. 
Using a density argument, we have that this holds for all $\phi\in H_0(\curl,D)$, thus
\[
\|\curl(a^\ep\curl u_1^\ep)+b^\ep u_1^\ep-\curl(a^0\curl u_0)-b^0u_0\|_{V'}\le c\ep
\]
so
\[
\|\curl(a^\ep\curl (u_1^\ep-u^\ep))+b^\ep (u_1^\ep-u^\ep)\|_{V'}\le c\ep.
\]
Let $\tau^\ep(x)$ be a function in ${\cal D}(D)$ such that $\tau^\ep(x)=1$ outside an $\ep$ neighbourhood of $\partial D$ and ${\rm sup}_{x\in D}\ep|\nabla\tau^\ep(x)|<c$ where $c$ is independent of $\ep$. Let
\[
w_1^\ep=u_0(x)+\ep\tau^\ep(x)N^r(x,{x\over\ep})\curl u_0(x)_r+\ep\nabla[\tau^\ep(x){\frak u}_1(x,{x\over\ep})].
\]
The function $w_1^\ep(x)$ belongs to $H_0(\scurl,D)$. We note that
\[
u_1^\ep-w_1^\ep=\ep(1-\tau^\ep(x))N^r(x,{x\over\ep})\curl u_0(x)_r+\ep\nabla[(1-\tau^\ep(x))\frak u_1(x,{x\over\ep})].
\]
From this,
\beqas
\curl(u_1^\ep-w_1^\ep)=\ep\curl_xN^r(x,{x\over\ep})\curl u_0(x)_r(1-\tau^\ep(x))+\curl_yN^r(x,{x\over\ep})\curl u_0(x)_r(1-\tau^\ep(x))-\\
\ep\curl u_0(x)_r\nabla\tau^\ep(x)\times N^r(x,{x\over\ep})+\ep(1-\tau^\ep(x))\nabla\curl u_0(x)_r\times N^r(x,{x\over\ep}).
\eeqas

Let $D^\ep\subset D$ be the $\ep$ neighbourhood of the boundary $\partial D$. We note that
\[
\|\phi\|_{L^2(D^\ep)}^2\le c\ep^2\|\phi\|_{H^1(D)}^2+c\ep\|\phi\|_{L^2(\partial D)}^2\le c\ep\|\phi\|_{H^1(D)}^2
\]
(see Hoang and Schwab \cite{HSmultirandom}).
We therefore deduce that $\|\curl u_0(x)_r\|_{L^2(D^\ep)}\le c\ep^{1/2}$. From these we have
\[
\|\curl(u_1^\ep-w_1^\ep)\|_{L^2(D)^3}\le c\ep^{1/2}.
\]
On the other hand,
\beqas
&& u_1^\ep-w_1^\ep=\ep(1-\tau^\ep(x))N^r(x,{x\over\ep})\curl u_0(x)_r-\ep\nabla\tau^\ep(x)w^r(x,{x\over\ep})u_{0r}(x)+\\
&& \ep(1-\tau^\ep(x))\nabla_xw^r(x,{x\over\ep})u_{0r}(x)+(1-\tau^\ep(x))\nabla_yw^r(x,{x\over\ep})u_{0r}(x)+\ep(1-\tau^\ep(x))w^r(x,{x\over\ep})\nabla u_{0r}(x).
\eeqas
Using the fact that $N^r\in C^1(\bar D,C(\bar Y))$, $w^r\in C^1(\bar D,C^1(\bar Y))$ and $\|u_{0r}\|_{L^2(D^\ep)^3}\le c \ep^{1/2}$, we deduce that $\|u_1^\ep-w_1^\ep\|_{L^2(D)^3}\le c\ep^{1/2}$. Therefore
\[
\|\curl(a^\ep\curl(u_1^\ep-w_1^\ep))+b^\ep(u_1^\ep-w_1^\ep)\|_{V'}\le c\ep^{1/2}.
\]
Thus
\[
\|\curl(a^\ep\curl(u^\ep-w_1^\ep))+b^\ep(u^\ep-w_1^\ep)\|_{V'}\le c\ep^{1/2}.
\]
Since $u^\ep-w_1^\ep\in H_0(\curl, D)$, $\|u^\ep-w_1^\ep\|_{H(\curl, D)}\le c\ep^{1/2}$. From these we have $\|u^\ep-u_1^\ep\|_{H(\curl,D)}\le c\ep^{1/2}$. We then get the conclusion. \eproof
\section{Homogenization error when $u_0\in H^s(\curl,D)$ for $0<s<1$}
In this section, we consider the case where the solution $u_0$ of the homogenized problem \eqref{eq:homprob} belongs to the weaker regularity space $H^s(\curl,D)$ for $0<s<1$. We have the following result.
\begin{theorem}\label{thm:main2}
Assume that $a\in C(\bar D,C(\bar Y))$, $u_0\in H^s(\curl,D)$, $N^r\in C^1(\bar D, C(\bar Y))^3$, $\scurl_yN^r\in C^1(\bar D,C(\bar Y))^3$, and $w^r\in C^1(\bar D,C^1(\bar Y))$ for all $r=1,2,3$, then 
\[
\|\ue-[u_0+\nabla_y\frak u_1(\cdot,{\cdot\over\ep})\|_{L^2(D)^3}\le c\ep^{s/(1+s)}
\]
and 
\[
\|\curl\ue-[\curl u_0+\scurl_yu_1(\cdot,{\cdot\over\ep})]\|_{L^2(D)^3}\le c\ep^{s/(1+s)}.
\]
\end{theorem}
\bproof
We consider a set of $M$ open cubes $Q_i$ ($i=1,\ldots,M$) of size $\ep^t$ for $t>0$ to be chosen later such that $D\subset\bigcup_{i=1}^MQ_i$ and $Q_i\bigcap D\ne\emptyset$. Each cube $Q_i$ intersects with only a finite number, which does not depend on $\ep$, of other cubes. We consider a partition of unity that consists of $M$ functions $\rho_i$ such that $\rho_i$ has support in $Q_i$, $\sum_{i=1}^M\rho_i(x)=1$ for all $x\in D$ and $|\nabla\rho_i(x)|\le c\ep^{-t}$ for all $x$ (indeed such a set of cubes $Q_i$ and a partition of unity can be constructed from a fixed set of cubes of size $O(1)$ by rescaling). For $r=1,2,3$ and $i=1,\ldots,M$, we denote by 
\[
U_i^r={1\over |Q_i|}\int_{Q_i}\curl u_0(x)_r dx
\]
and
\[
V_i^r={1\over |Q_i|}\int_{Q_i}u_0(x)_rdx
\]
(as $u_0\in H^s(D)^3$ and $\curl u_0\in H^s(D)^3$, for the Lipschitz domain $D$, we can extend each of them, separately, continuously outside $D$ and understand $u_0$ and $\curl u_0$ as these extensions (see Wloka \cite{Wloka} Theorem 5.6)).
Let $U_i$ and $V_i$ denote the vector $(U_i^1, U_i^2,U_i^3)$ and $(V_i^1,V_i^2,V_i^3)$ respectively. 
Let $B$ be the unit cube in $\IR^3$. From Poincare inequality, we have
\[
\int_B\Big|\phi-\int_B\phi(x)dx\Big|^2dx\le c\int_B|\nabla \phi(x)|^2dx\ \ \forall\,\phi\in H^1(B).
\]
By translation and scaling, we deduce that 
\[
\int_{Q_i}\Big|\phi-{1\over |Q_i|}\int_{Q_i}\phi(x)dx\Big|^2dx\le c\ep^{2t}\int_{Q_i}|\nabla\phi(x)|^2dx\ \ \forall\,\phi\in H^1(Q_i)
\]
i.e
\[
\Big\|\phi-{1\over |Q_i|}\int_{Q_i}\phi(x)dx\Big\|_{L^2(Q_i)}\le c\ep^t\|\phi\|_{H^1(Q_i)}.
\]
Together with
\[
\Big\|\phi-{1\over |Q_i|}\int_{Q_i}\phi(x)dx\Big\|_{L^2(Q_i)}\le c\|\phi\|_{L^2(Q_i)}
\]
we deduce from interpolation that
\[
\Big\|\phi-{1\over |Q_i|}\int_{Q_i}\phi(x)dx\Big\|_{L^2(Q_i)}\le c\ep^{ts}\|\phi\|_{H^s(Q_i)}\ \ \forall\,\phi\in H^s(Q_i).
\]
Thus
\be
\int_{Q_i}|\curl u_0(x)_r-U_i^r|^2 dx\le c\ep^{2ts}\|(\curl u_0)_r\|^2_{H^s(Q_i)}.
\label{eq:intine}
\ee
Let 
\[
u_1^\ep(x)=u_0(x)+\ep N^r(x,{x\over \ep})U_j^r\rho_j(x)+\ep\nabla [w^r(x,{x\over\ep})V_j^r\rho_j(x)].
\]
We have 
\beqas
&&\curl(a^\ep(x)\curl u_1^\ep(x))+b^\ep(x)u_1^\ep(x)\\
&&\qquad=\curl a(x,{x\over \ep})\left[\curl u_0(x)+\ep\scurl_xN^r(x,{x\over\ep})U_j^r\rho_j(x)+\scurl_yN^r(x,{x\over\ep})U_j^r\rho_j+\right.\\
&&\qquad\qquad\qquad\left.\ep(U_j^r\nabla\rho_j)\times N^r(x,{x\over\ep})\right]+\\
&&\qquad\qquad b(x,{x\over\ep})\left[u_0(x)+\ep N^r(x,{x\over\ep})U_j^r\rho_j(x)+\ep\nabla_x w^r(x,{x\over\ep})V_j^r\rho_j(x)+\nabla_yw^r(x,{x\over\ep})V_j^r\rho_j(x)+\right.\\
&&\qquad\qquad\qquad\left.\ep w^r(x,{x\over\ep})V_j^r\nabla\rho_j(x)\right]\\
&&=\curl(a^0(x)\curl u_0(x))+b^0(x)u_0(x)+\scurl\left[G_r(x,{x\over\ep})U_j^r\rho_j(x)\right]+g_r(x,{x\over\ep})V^r_j\rho_j(x)+\ep\curl I(x)+\ep J(x)+\\
&&\qquad\qquad\qquad \curl[(a^\ep(x)-a^0(x))(\curl u_0(x)-U_j\rho_j(x))]+(b^\ep(x)-b^0(x))(u_0(x)-V_j\rho_j(x)),
\eeqas
where $G_r(x,y)$ and $g_r(x,y)$ are defined as in \eqref{eq:G} and \eqref{eq:g} respectively, and 
\[
I(x)=a(x,{x\over \ep})\left[\scurl_xN^r(x,{x\over\ep})U_j^r\rho_j(x)+(U_j^r\nabla\rho_j(x))\times N^r(x,{x\over\ep})\right],
\]
and 
\[
J(x)=b(x,{x\over\ep})\left[N^r(x,{x\over\ep})U_j^r\rho_j(x)+\nabla_xw^r(x,{x\over\ep})V_j^r\rho_j(x)+w^r(x,{x\over\ep})V_j^r\nabla\rho_j(x)\right].
\]
Therefore for $\phi\in H_0(\curl,D)$
\beqas
&&\langle\curl(a^\ep\curl u_1^\ep)+b^\ep u_1^\ep -\curl(a^0\curl u_0)-b^0u_0,\phi\rangle\\
&&=\int_DG_r(x,{x\over\ep})U_j^r\rho_j(x)\cdot\curl\phi dx+\int_Dg_r(x,{x\over\ep})(V_j^r\rho_j(x))\cdot\phi(x)dx+\\
&&\ep\int_DI(x)\cdot\curl\phi(x)dx+\ep\int_DJ(x)\cdot\phi(x)dx+\int_D(a^\ep-a^0)(\curl u_0(x)-U_j\rho_j)\cdot\curl\phi(x)dx\\
&&\qquad\qquad+\int_D(b^\ep-b^0)(u_0-V_j\rho_j)\cdot\phi dx.
\eeqas
We have that
\beqas
&&\int_DG_r(x,{x\over\ep})(U_j^r\rho_j)\cdot\curl\phi dx=\int_D\left[\ep\nabla\tilde G_r(x,{x\over\ep})-\ep\nabla_x\tilde G_r(x,{x\over\ep})\right](U^r_j\rho_j(x))\cdot\curl\phi dx\\
&&\qquad\qquad=-\ep\int_D\tilde G_r(x,{x\over\ep}){\rm div}[(U_j^r\rho_j)\curl\phi]dx-\ep\int_D\nabla_x\tilde G_r(x,{x\over\ep})(U_j^r\rho_j)\cdot\curl\phi dx.
\eeqas
We note that
\beqas
\Big|\int_D\nabla_x\tilde G_r(x,{x\over\ep})(U_j^r\rho_j)\cdot\curl\phi dx\Big|\le c\|(U_j^r\rho_j)\|_{L^2(D)}\|\curl\phi\|_{L^2(D)^3}.
\eeqas
From 
\[
\|U_j^r\rho_j\|_{L^2(D)}^2=\int_D(U_j^r)^2\rho_j(x)^2dx+\sum_{i\ne j}\int_D U_i^rU_j^r\rho_i(x)\rho_j(x)dx,
\]
and the fact that the support of each function $\rho_i$ intersects only with the support of a finite number (which does not depend on $\ep$) of other functions $\rho_j$ in the partition of unity, we deduce 
\beqas
\|U_j^r\rho_j\|^2_{L^2(D)}\le c\sum_{j=1}^M(U_j^r)^2|Q_j|=c\sum_{j=1}^M{1\over|Q_j|}\left(\int_{Q_j}\curl u_0(x)_rdx\right)^2\le c\sum_{j=1}^M\int_{Q_j}\curl u_0(x)_r^2dx\le c\int_D\curl u_0(x)_r^2dx.
\eeqas
Thus 
\beqas
\Big|\ep\int_D\nabla_x\tilde G_r(x,{x\over\ep})(U_j^r\rho_j)\cdot\curl\phi dx\Big|\le c\ep\|\curl\phi\|_{L^2(D)^3}.
\eeqas
We also have
\beqas
\ep\int_D\tilde G_r(x,{x\over\ep}){\rm div}[(U_j^r\rho_j)\curl\phi]dx=\ep\int_D\tilde G_r(x,{x\over\ep})\left[(U_j^r\nabla\rho_j(x))\right]\cdot\curl\phi dx\le c\ep\|U_j^r\nabla\rho_j\|_{L^2(D)^3}\|\curl\phi\|_{L^2(D)^3}.
\eeqas
As the support of each function $\rho_i$ intersects with the support of a finite number of other functions $\rho_j$ and $\|\nabla\rho_j\|_{L^\infty(D)}\le c\ep^{-t}$, we have
\[
\|U_j^r\nabla\rho_j\|_{L^2(D)^3}^2\le c\sum_{j=1}^M(U_j^r)^2|Q_j|\|\nabla\rho_j\|^2_{L^\infty(D)}\le c\ep^{-2t}\sum_{j=1}^M(U^r_j)^2|Q_j|\le c\ep^{-2t},
\]
so
\[
\ep\int_D\tilde G_r(x,{x\over\ep}){\rm div}[(U_j^r\rho_j)\curl\phi]dx\le c\ep\|U_j^r\nabla\rho_j\|_{L^2(D)^3}\|\curl\phi\|_{L^2(D)^3}\le c\ep^{1-t}\|\curl\phi\|_{L^2(D)^3}.
\]
We therefore deduce that
\[
\Big|\int_DG_r(x,{x\over\ep})(U_j^r\rho_j)\cdot\curl\phi dx\Big|\le c\ep^{1-t}\|\curl\phi\|_{L^2(D)^3}.
\]
We have
\[
\int_Dg_r(x,{x\over\ep})(V_j^r\rho_j)\cdot\phi(x)dx=\int_D[\ep\curl\tilde g_r(x,{x\over\ep})-\ep\scurl_x\tilde g_r(x,{x\over\ep})](V_j^r\rho_j)\cdot\phi dx.
\]
Arguing similarly as above, we have
\[
\Big|\ep\int_D\scurl_x\tilde g_r(x,{x\over\ep}))(V_j^r\rho_j)\cdot\phi dx\Big|\le c\ep\|V_j^r\rho_j\|_{L^2(D)^3}\|\phi\|_{L^2(D)^3}\le c\ep\|\phi\|_{L^2(D)^3},
\]
and 
\beqas
\Big|\ep\int_D\curl\tilde g_r(x,{x\over\ep})(V_j^r\rho_j)\cdot\phi dx\Big|&\le& \Big|\ep\int_D\tilde g_r(x,{x\over\ep})\cdot\curl[(V_j^r\rho_j)\phi] dx\Big|\\
&\le& \Big|\ep\int_D\tilde g_r(x,{x\over\ep})\cdot[(V_j^r\rho_j)\curl\phi+\phi\times(V_j^r\nabla\rho_j)]dx\\
&\le& c\left(\ep\|\curl\phi\|_{L^2(D)^3}+c\ep^{1-t}\|\phi\|_{L^2(D)^3}\right)\left(\sum_{j=1}^M(V_j^r)^2|Q_j|\right)^{1/2}\\
&\le& c(\ep\|\curl\phi\|_{L^2(D)^3}+c\ep^{1-t}\|\phi\|_{L^2(D)^3}).
\eeqas
We note that
\[
\|I\|_{L^2(D)^3}\le c\sup_r[\|U_j^r\rho_j\|_{L^2(D)}+\|U_j^r\nabla\rho_j\|_{L^2(D)}]\le c\ep^{-t},
\]
and 
\[
\|J\|_{L^2(D)^3}\le c\sup_r[\|U_j^r\rho_j\|_{L^2(D)}+c\|V_j^r\rho_j\|_{L^2(D)}+c\|V_j^r\nabla\rho_j\|_{L^2(D)}]\le c\ep^{-t}.
\]
We have further that
\[
\langle\curl((a^\ep-a^0)(\curl u_0-U_j\rho_j)),\phi\rangle\le c\|\curl u_0-(U_j\rho_j))\|_{L^2(D)^3}\|\curl\phi\|_{L^2(D)^3}.
\]
We note that
\beqas
\int_D|(\curl u_0)_r-(U_j^r\rho_j)|^2dx=\int_D|\sum_{j=1}^M((\curl u_0)_r-U_j^r)\rho_j|^2dx.
\eeqas
Using the support property of $\rho_j$, we have from \eqref{eq:intine}
\begin{eqnarray}
\int_D|(\curl u_0)_r-(U_j^r\rho_j)|^2dx&\le& c\sum_{j=1}^M\int_{Q_j}|(\curl u_0)_r-U_j^r|^2dx\le c\ep^{2st}\sum_{j=1}^M\|(\curl u_0)_r\|_{H^s(Q_j)}^2\nonumber\\
&=&c\ep^{2st}\sum_{j=1}^M\left[\int_{Q_j}(\curl u_0)_r^2dx+\int_{Q_j\times Q_j}{(\curl u_0(x)_r-\curl u_0(x')_r)^2\over |x-x'|^{2s}}dxdx'\right]\nonumber\\
&\le& c\ep^{2st}\left[\|(\curl u_0)_r\|_{L^2(D)}^2+\int_{D\times D}{(\curl u_0(x)_r-\curl u_0(x')_r)^2\over |x-x'|^{2s}}dxdx'\right]\nonumber\\
&=&c\ep^{2st}\|(\curl u_0)_r\|_{H^s(D)}^2.
\label{eq:2}
\end{eqnarray}
Thus 
\[
\langle\curl((a^\ep-a^0)(\curl u_0-U_j\rho_j)),\phi\rangle\le c\ep^{st}\|\curl\phi\|_{L^2(D)^3}.
\]
Similarly, we have
\[
\Big|\int_D(b^\ep-b^0)(u_{0}-\sum_{j=1}^MV_j\rho_j)\cdot\phi dx\Big|\le c\|\sum_{j=1}^M(u_{0}-V_j)\rho_j\|_{L^2(D)^3}\|\phi\|_{L^2(D)^3}\le c\ep^{st}\|\phi\|_{L^2(D)^3}.
\]
Therefore
\[
|\langle\curl(a^\ep\curl u_1^\ep)+b^\ep u_1^\ep-\curl(a^0\curl u_0)-b^0u_0,\phi\rangle|\le c(\ep^{1-t}+\ep^{st})\|\phi\|_V
\]
i.e.
\[
\|\curl(a^\ep\curl u_1^\ep)+b^\ep u_1^\ep-\curl(a^0\curl u_0)-b^0u_0\|_{V'}\le c(\ep^{1-t}+\ep^{st}).
\]
Thus
\be
\|\curl(a^\ep\curl u_1^\ep)+b^\ep u_1^\ep-\curl(a^\ep\curl u^\ep)-b^\ep u^\ep\|_{V'}\le c(\ep^{1-t}+\ep^{st}).
\label{eq:1}
\ee
We choose $\tau^\ep$ as in the previous section and consider the function
\[
w_1^\ep(x)=u_0(x)+\ep\tau^\ep(x)N^r(x,{x\over\ep})U_j^r\rho_j(x)+\ep\nabla[\tau^\ep(x)w^r(x,{x\over\ep})V_j^r\rho_j].
\]
We then have
\[
u_1^\ep-w_1^\ep=\ep(1-\tau^\ep(x))N^r(x,{x\over\ep})U_j^r\rho_j(x)+\ep\nabla[(1-\tau^\ep(x))w^r(x,{x\over\ep})V_j^r\rho_j]
\]
and
\beqas
\curl(u_1^\ep-w_1^\ep)=\ep\scurl_xN^r(x,{x\over\ep})U_j^r\rho_j(x)(1-\tau^\ep(x))+\scurl_yN^r(x,{x\over\ep})U_j^r\rho_j(x)(1-\tau^\ep(x))-\\
\ep U_j^r\rho_j(x)\nabla\tau^\ep(x)\times N^r(x,{x\over\ep})+\ep(1-\tau^\ep(x))U_j^r\nabla\rho_j(x)\times N^r(x,{x\over\ep}).
\eeqas
As shown above $\|U_j^r\rho_j\|_{L^2(D)}$ is uniformly bounded, so 
\[
\|\ep\scurl_xN^r(x,{x\over\ep})(U_j^r\rho_j)(1-\tau^\ep(x))\|_{L^2(D)^3}\le c\ep.
\]
Let $\tilde D^\ep$ be the $3\ep^{t}$ neighbourhood of $\partial D$. We note that $\curl u_0$ is extended continuously outside $D$. 
As shown in Hoang and Schwab \cite{HSmultirandom}, for $\phi\in H^1(\tilde D^\ep)$
\[
\|\phi\|_{L^2(\tilde D^\ep)}\le c\ep^{t/2}\|\phi\|_{H^1(\tilde D^\ep)}.
\]
From this and 
\[
\|\phi\|_{L^2(\tilde D^\ep)}\le \|\phi\|_{L^2(\tilde D^\ep)},
\]
using interpolation we get 
\[
\|\phi\|_{L^2(\tilde D^\ep)}\le c\ep^{st/2}\|\phi\|_{H^s(\tilde D^\ep)}\le c\ep^{st/2}\|\phi\|_{H^s(D)},
\]
$\forall\,\phi\in H^s(D)$ extended continuously outside $D$.

We then have
\beqas
\|U_j^r\rho_j\|^2_{L^2(D^\ep)}&\le& c\sum_{j=1}^M\int_{Q_j\bigcap D^\ep}(U_j^r)^2\rho_j^2dx\\
&\le& c\sum_{j=1}^M|Q_j\bigcap D^\ep|{1\over |Q_j|^2}(\int_{Q_j}(\curl u_0)_rdx)^2\\
&\le& c\sum_{Q_j\bigcap D^\ep\ne\emptyset}{|Q_j\bigcap D^\ep|\over|Q_j|}\int_{Q_j}(\curl u_0)_r^2dx.
\eeqas
As $D^\ep$ is the $\ep$ neighbourhood of $\partial D$ and $Q_j$ has size $\ep^t$, $|Q_j\bigcap D^\ep|\le c\ep^{1+(d-1)t}$ so $|Q_j\bigcap D^\ep|/|Q_j|\le c\ep^{1-t}$. When $Q_j\bigcap D^\ep\ne\emptyset$, $Q_j\subset\tilde D^\ep$. Thus
\beqas
\|U_j^r\rho_j\|^2_{L^2(D^\ep)}\le
c\ep^{1-t}\|(\curl u_0)_r\|^2_{L^2(\tilde D^\ep)}\le c\ep^{1-t+st}\|\curl u_0\|^2_{H^s(D)^3}.
\eeqas 
Therefore
\[
\|\scurl_yN^r(x,{x\over\ep})(U_j^r\rho_j)(1-\tau^\ep(x))\|_{L^2(D)^3}\le c\ep^{(1-t+st)/2}
\]
and
\[
\|\ep(U_j^r\rho_j)\nabla\tau^\ep(x)\times N^r(x,{x\over\ep})\|_{L^2(D)^3}\le c\ep^{(1-t+st)/2}.
\]
Similarly we have
\beqas
\|U_j^r\nabla\rho_j\|_{L^2(D^\ep)^3}^2\le c\ep^{-2t}\sum_{Q_j\bigcap D^\ep\ne\emptyset}|Q_j\bigcap D^\ep|(U_j^r)^2\le  c\ep^{-2t}\sum_{Q_j\bigcap D^\ep\ne\emptyset}{|Q_j\bigcap D^\ep|\over|Q_j|}\int_{Q_j}(\curl u_0)_r^2dx\\
\le c\ep^{-2t+1-t}\|\curl u_0\|_{L^2(\tilde D^\ep)^3}^2\le c\ep^{1-3t+st}\|\curl u_0\|^2_{H^s(D)^3}.
\eeqas
Thus 
\[
\|\ep(1-\tau^\ep(x))(U_j^r\nabla\rho_j)\times N^r(x,{x\over\ep})\|_{L^2(D)}\le c\ep^{(1-t)+(1-t+st)/2}.
\]
Therefore
\[
\|\curl(u_1^\ep-w_1^\ep)\|_{L^2(D)^3}\le c(\ep^{(1-t+st)/2}+\ep^{(1-t)+(1-t+st)/2}).
\]
We further have that
\beqas
\ep\nabla[(1-\tau^\ep(x))w^r(x,{x\over\ep})(V_j^r\rho_j)]=-\ep\nabla\tau^\ep(x)w^r(x,{x\over\ep})(V_j^r\rho_j)+\ep(1-\tau^\ep(x))\nabla_xw^r(x,{x\over\ep})(V_j^r\rho_j)\\
+(1-\tau^\ep(x))\nabla_yw^r(x,{x\over\ep})(V_j^r\rho_j)+\ep(1-\tau^\ep(x))w^r(x,{x\over\ep})(V_j^r\nabla\rho_j).
\eeqas
Arguing as above, we deduce that
\[
\|V_j^r\rho_j\|_{L^2(D^\ep)}\le c\ep^{(1-t+st)/2},\ \ \|V_j^r\nabla\rho_j\|_{L^2(D^\ep)}\le c\ep^{(1-t+st)/2-t}.
\]
Therefore
\[
\|\ep\nabla[(1-\tau^\ep(x))w^r(x,{x\over\ep})(V_j^r\rho_j)]\|_{L^2(D)^3}\le c(\ep^{(1-t+st)/2}+\ep^{1-t+(1-t+st)/2}).
\]
Thus
\[
\|u_1^\ep-w_1^\ep\|_{L^2(D)^3}\le c(\ep^{(1-t+st)/2}+\ep^{(1-t)+(1-t+st)/2}).
\]
Choosing $t=1/(s+1)$ we have
\[
\|\curl(a^\ep\curl(u_1^\ep-w_1^\ep))+b^\ep(u_1^\ep-w_1^\ep)\|_{V'}\le c\ep^{s/(s+1)}.
\]
This together with \eqref{eq:1} gives
\[
\|\curl(a^\ep\curl(u^\ep-w_1^\ep))+b^\ep(u^\ep-w_1^\ep)\|_{V'}\le c\ep^{s/(s+1)}.
\]
Thus
\[
\|\ue-w_1^\ep\|_{V}\le c\ep^{s/(s+1)}
\]
which implies
\be
\|\ue-u_1^\ep\|_V\le c\ep^{s/(s+1)}.
\label{eq:eqs}
\ee
We note that 
\[
\curl u_1^\ep=\curl u_0(x)+\scurl_yN^r(x,{x\over\ep})(U_j^r\rho_j)+\ep\scurl_xN^r(x,{x\over\ep})(U_j^r\rho_j)+\ep N^r(x,{x\over\ep})\times (U_j^r\nabla\rho_j).
\]
From
\[
\|\ep\scurl_xN^r(x,{x\over\ep})(U_j^r\rho_j)\|_{L^2(D)^3}\le c\ep,\ \ \mbox{and}\ \ \|\ep N^r(x,{x\over\ep})\times (U_j^r\nabla\rho_j)\|_{L^2(D)^3}\le c\ep\ep^{-t}=c\ep^{s/(1+s)},
\]
we deduce that
\[
\|\curl\ue_1-\curl u_0-\scurl_yN^r(x,{x\over\ep})(U_j^r\rho_j)\|_{L^2(D)^3}\le c\ep^{s/(s+1)}.
\]
From \eqref{eq:2}
\[
\|\curl u_0-(U_j\rho_j)\|_{L^2(D)^3}\le c\ep^{ts}= c\ep^{s/(s+1)},
\]
we get
\[
\|\curl u^\ep_1-[\curl u_0+\scurl_yN^r(x,{x\over\ep})(\curl u_0)_r\|_{L^2(D)^3}\le c\ep^{s/(s+1)}.
\]
This together with \eqref{eq:eqs} implies
\[
\|\curl u^\ep-[\curl u_0+\scurl_yN^r(x,{x\over\ep})(\curl u_0)_r\|_{L^2(D)^3}\le c\ep^{s/(s+1)}.
\]
\eproof

\section{Regularity of $N^r$, $w^r$ and $u_0$}
We show in this section that the regularity requirements for the solutions $N^r(x,y)$, $w^r(x,y)$ of the cell problems \eqref{eq:cellprobN} and \eqref{eq:cellprobw}, and of the solution $u_0$ of the homogenized problem \eqref{eq:homprob} hold. We first prove the following lemma.
\begin{lemma}\label{lem:lem}
Let $\psi\in H_\#(\curl,Y)\bigcap H_\#({\rm div},Y)$. Assume further that $\int_Y\psi(y)dy=0$. Then $\psi\in H^1_\#(Y)$ and 
\[
\|\psi\|_{H^1(Y)^3}\le c(\|\scurl_y\psi\|_{L^2(Y)^3}+\|{\rm div}\psi\|_{L^2(Y)^3}).
\]
\end{lemma}
\bproof
Let $\omega\subset\IR^3$ be a smooth domain such that $\omega\supset Y$. Let $\eta\in {\cal D}(\omega)$ be such that $\eta(y)=1$ when $y\in Y$. We have 
\[
\scurl_y(\eta\psi)=\eta\scurl_y\psi+\nabla_y\eta\times\psi\in L^2(\omega)^3
\]
 and 
\[
{\rm div}_y(\eta\psi)=\nabla_y\eta\cdot\psi+\eta{\rm div}_y\psi\in L^2(\omega)^3.
\]
Together with the boundary condition, we conclude that $\eta\psi\in H^1(\omega)^3$ so $\psi\in H^1(Y)^3$.

We note that
\beqas
\int_Y({\rm div}\psi(y)^2+|\scurl_y\psi(y)|^2)dy=\sum_{i,j=1}^3\int_Y\left({\partial\psi_i\over\partial y_j}\right)^2+\sum_{i\ne j}\int_Y{\partial\psi_i\over\partial y_i}{\partial\psi_j\over\partial y_j}dy-\sum_{i\ne j}\int_Y{\partial\psi_j\over\partial y_i}{\partial\psi_i\over\partial y_j}dy.
\eeqas
Assume that $\psi$ is a smooth periodic function. We have
\beqas
\int_Y{\partial\psi_i\over\partial y_i}{\partial\psi_j\over\partial y_j}dy=\int_Y\left[{\partial\over\partial y_i}\left(\psi_i{\partial\psi_j\over\partial y_j}\right)-\psi_i{\partial^2\psi_j\over\partial y_i\partial y_j}\right]dy=-\int_Y\psi_i{\partial^2\psi_j\over\partial y_i\partial y_j}dy
\eeqas
as $\psi$ is periodic. Similarly, we have
\beqas
\int_Y{\partial\psi_i\over\partial y_j}{\partial\psi_j\over\partial y_i}dy=\int_Y\left[{\partial\over\partial y_j}\left(\psi_i{\partial\psi_j\over\partial y_i}\right)-\psi_i{\partial^2\psi_j\over\partial y_j\partial y_i}\right]dy=-\int_Y\psi_i{\partial^2\psi_j\over\partial y_j\partial y_i}dy.
\eeqas
Thus 
\[
\int_Y{\partial\psi_i\over\partial y_i}{\partial\psi_j\over\partial y_j}dy=\int_Y{\partial\psi_i\over\partial y_j}{\partial\psi_j\over\partial y_i}dy.
\]
Therefore
\[
\|\nabla_y\psi\|^2_{L^2(Y)^3}= \|{\rm div}\psi\|^2_{L^2(Y)}+\|\scurl_y\psi\|^2_{L^2(Y)^3}.
\]
Using a density argument, this holds for all $\phi\in H^1_\#(Y)$.
As $\int_Y\psi(y)dy=0$, from Poincare inequality we deduce
\[
\|\psi\|_{H^1(Y)}\le c(\|{\rm div}\psi\|_{L^2(Y)}+\|\scurl_y\psi\|_{L^2(Y)^3}).
\]
\eproof

We then prove the following result on the regularity of $N^r$.

\begin{proposition}\label{prop:5.2}
Assume that $a(x,y)\in C^1(\bar D,C^2(\bar Y))^{3\times 3}$ and is isotropic, then $\scurl_yN^r(x,y)\in C^1(\bar D,C(\bar Y))^3$ and we can choose a version of $N^r$ in $L^2(D,\tilde H_\#(\curl,Y))$ so that $N^r(x,y)\in C^1(\bar D, C(\bar Y))^3$. 
\end{proposition}
\bproof
As $a$ is isotropic, we understand in this proof that it is a scalar function. 
From \eqref{eq:cellprobN}, we have 
\[
\scurl_y(\scurl_y N^r)(x,\cdot)=-{1\over a(x,\cdot)}\nabla_y a(x,\cdot)\times\scurl_yN^r(x,\cdot)-{1\over a(x,\cdot)}\scurl_y(a(x,\cdot)e_r)\in L^2(Y)^3.
\]
Therefore $\scurl_yN^r\in H_\#(\curl,Y)$. This together with ${\rm div}_y\scurl_yN^r=0$ and Lemma \ref{lem:lem} implies $\scurl_yN^r\in H^1_\#(Y)^3$. 

Fixing an index $q=1,2,3$, we have
\[
\scurl_y{\partial\over\partial y_q}\scurl_yN^r={\partial\over\partial y_q}\scurl_y\scurl_y N^r=-{\partial\over\partial y_q}\left({1\over a}\nabla_ya\times\scurl_yN^r\right)-{\partial\over\partial y_q}\left({1\over a}\scurl_y(ae_r)\right)\in L^2(Y)^3.
\]
With ${\rm div}_y{\partial\over\partial y_q}\scurl_yN^r=0$, we deduce that
${\partial\over\partial y_q}\scurl_yN^r(x,\cdot)\in H^1(Y)^3$ for all $q=1,2,3$ so $\scurl_yN^r(x,\cdot)\in H^2(Y)^3\subset C(\bar Y)^3$. 

We now show that $\scurl_yN^r\in C^1(\bar D,H^2(Y))^3\subset C^1(\bar D,C(\bar Y))^3$. 

Fix $h\in\IR^3$. From \eqref{eq:cellprobN} we have
\beqas
\scurl_y(a(x,y)\scurl_y(N^r(x+h,y)-N^r(x,y)))=-\scurl_y((a(x+h,y)-a(x,y))e_r)-\\
\scurl_y((a(x+h,y)-a(x,y))\scurl_yN^r(x+h,y)).
\eeqas
This together with 
the smoothness of $a$ and the uniformly boundedness of $\scurl_yN^r(x,y)$ in $L^2(Y)^3$ gives
\be
\lim_{h\to 0}\|\scurl_y(N^r(x+h,\cdot)-N^r(x,\cdot))\|_{L^2(Y)^3}=0.
\label{eq:e1}
\ee
We also have that
\begin{eqnarray}
\scurl_y\scurl_y(N^r(x+h,y)-N^r(x,y))=-{1\over a}\nabla_y a\times \scurl_y((N^r(x+h,y)-N^r(x,y)))\nonumber\\
-{1\over a}\scurl_y((a(x+h,y)-a(x,y))e_r)\nonumber
\\
-{1\over a}\scurl_y((a(x+h,y)-a(x,y))\scurl_yN^r(x+h,y))
\label{eq:e11}
\end{eqnarray}
so
\be
\lim_{|h|\to 0}\|\scurl_y\scurl_y(N^r(x+h,\cdot)-N^r(x,\cdot))\|_{L^2(Y)^3}=0.
\label{eq:e2}
\ee
From Lemma \ref{lem:lem}, we have
\be
\lim_{|h|\to 0}\|\scurl_y(N^r(x+h,\cdot)-N^r(x,\cdot))\|_{H^1(Y)^3}=0.
\label{eq:e21}
\ee
From \eqref{eq:e11}, we get
\be
\lim_{|h|\to 0}\|\scurl_y\scurl_y(N^r(x+h,\cdot)-N^r(x,\cdot))\|_{H^1(Y)^3}=0.
\label{eq:e22}
\ee
We have
\beqas
\scurl_y\left[a(x,y)\scurl_y({N^r(x+h,y)-N^r(x,y)\over h})\right]=-\scurl_y(({a(x+h,y)-a(x,y)\over h})e_r)-\\
\scurl_y\left({a(x+h,y)-a(x,y)\over h}\scurl_yN^r(x+h,y)\right).
\eeqas
Let $\chi^r(x,\cdot)\in \tilde H_\#(\curl,Y)$ be the solution of the problem
\beqas
\scurl_y(a(x,y)\scurl_y\chi^r(x,\cdot))=-\scurl_y({\partial a\over\partial x_q}e_r)-\scurl_y\left({\partial a\over\partial x_q}\scurl_y N^r(x,y)\right).
\eeqas
We deduce that
\begin{eqnarray}
&&\scurl_y\left(a(x,y)\scurl_y\left({N^r(x+h,y)-N^r(x,y)\over h}-\chi^r(x,y)\right)\right)\nonumber\\
&&\qquad=-\scurl_y\left(\left( {a(x+h,y)-a(x,y)\over h}-{\partial a\over\partial x_q}(x,y)\right)e_r\right)\nonumber\\
&&\qquad-\scurl_y\left(\left({a(x+h,y)-a(x,y)\over h}-{\partial a\over\partial x_q}\right)\scurl_yN^r(x+h,y)\right)\nonumber\\
&&\qquad-\scurl_y\left({\partial a\over\partial x_q}(\scurl_y(N^r(x+h,y)-N^r(x,y)))\right).
\label{eq:e3}
\end{eqnarray}
Let $h\in \IR^3$ be a vector whose all components are 0 except the $q$th component. 
We have 
\begin{eqnarray}
\Big\|\scurl_y\left({N^r(x+h,\cdot)-N^r(x,\cdot)\over h}-\chi^r(x,\cdot)\right)\Big\|_{L^2(Y)^3}\le
c\sup_{y\in Y}\Big\|{a(x+h,y)-a(x,y)\over h}-{\partial a\over\partial x_q}(x,y)\Big\|_{L^\infty(Y)}+\nonumber\\
c\|\scurl_y(N^r(x+h,\cdot)-N^r(x,\cdot)\|_{L^2(Y)^3}.
\label{eq:ee}
\end{eqnarray}
Thus from \eqref{eq:e1} we have
\be
\lim_{|h|\to 0}\Big\|\scurl_y\left({N^r(x+h,\cdot)-N^r(x,\cdot)\over h}-\chi^r(x,\cdot)\right)\Big\|_{L^2(Y)^3}=0.
\label{eq:e4}
\ee
From \eqref{eq:e3}, we have
\beqas
&&\scurl_y\scurl_y\left({N^r(x+h,y)-N^r(x,y)\over h}-\chi^r(x,y)\right)\\
&&\qquad=-{1\over a}\nabla_y a\times\scurl_y\left({N^r(x+h,y)-N^r(x,y)\over h}-\chi^r(x,y)\right)\\
&&\qquad-{1\over a}\scurl_y\left(\left({a(x+h,y)-a(x,y)\over h}-{\partial a(x,y)\over\partial x_q}\right)e_r\right)\\
&&\qquad-{1\over a}\scurl_y\left(\left({a(x+h,y)-a(x,y)\over h}-{\partial a\over\partial x_q}\right)\scurl_yN^r(x+h,y)\right)\\
&&\qquad-{1\over a}\scurl_y\left({\partial a\over\partial x_q}\scurl_y(N^r(x+h,y)-N^r(x,y))\right)
\eeqas
which converges to 0 in the $L^2(Y)$ norm when $|h|\to 0$ due to \eqref{eq:e1}, \eqref{eq:e21} and \eqref{eq:e4}. 
From Lemma \ref{lem:lem}, we have
\be
\lim_{|h|\to 0}\Big\|\scurl_y\left({N^r(x+h,\cdot)-N^r(x,\cdot)\over h}-\chi^r(x,\cdot)\right)\Big\|_{H^1(Y)^3}=0.
\label{eq:e5}
\ee
Let $p=1,2,3$. We then have
\beqas
&&{\partial\over\partial y_p}\scurl_y\scurl_y\left({N^r(x+h,y)-N^r(x,y)\over h}-\chi^r(x,y)\right)\\
&&\qquad=-{\partial\over\partial y_p}\left[{1\over a}\nabla a\times\scurl_y\left({N^r(x+h,y)-N^r(x,y)\over h}-\chi^r(x,y)\right)\right]\\
&&\qquad-{\partial\over\partial y_p}\left[{1\over a}\scurl_y\left(\left({a(x+h,y)-a(x,y)\over h}-{\partial a(x,y)\over\partial x_q}\right)e_r\right)\right]\\
&&\qquad-{\partial\over\partial y_p}\left[{1\over a}\scurl_y\left(\left({a(x+h,y)-a(x,y)\over h}-{\partial a\over\partial x_q}\right)\scurl_yN^r(x+h,y)\right)\right]\\
&&\qquad-{\partial\over\partial y_p}\left[{1\over a}\scurl_y\left({\partial a\over\partial x_q}\scurl_y(N^r(x+h,y)-N^r(x,y))\right)\right]
\eeqas
which converges to 0 in $L^2(Y)$ for each $x$ due to \eqref{eq:e22}, \eqref{eq:e5} and the uniform boundedness of $\|\curl N^r(x,\cdot)\|_{H^2(Y)^3}$. 
Hence
\[
\lim_{|h|\to 0}\Big\|{\partial\over\partial y_p}\scurl_y\left({N^r(x+h,\cdot)-N^r(x,\cdot)\over h}-\chi^r(x,\cdot)\right)\Big\|_{H^1(Y)^3}=0
\]
i.e.
\[
{\partial\over\partial x_q}\scurl_yN^r(x,\cdot)=\chi^r(x,\cdot)
\]
in the $H^2(Y)^3$ norm. Thus $\scurl_yN^r\in C^1(\bar D,H^2(Y))^3$.

We can always choose a version of $N^r$ such that ${\rm div}_yN^r=0$. Indeed, let $\Phi(x,\cdot)\in L^2(D, H^1_\#(Y))$ be a function such that $\Delta_y\Phi=-{\rm div}_yN^r$, then $\curl_y(N^r+\nabla_y\Phi)=\scurl_yN^r$ and ${\rm div}_y(N^r+\nabla_y\Phi)=0$. Thus we assume that $N^r$ and $\chi^r$ are chosen so that they are divergence free with respect to $y$. Further, we choose them so that $\int_YN^r(x,y)dy=0$ and $\int_Y\chi^r(x,y)dy=0$. From \eqref{eq:e5} and Lemma \ref{lem:lem}, we have
\[
\lim_{|h|\to 0}\Big\|\left({N^r(x+h,\cdot)-N^r(x,\cdot)\over h}-\chi^r(x,\cdot)\right)\Big\|_{H^1(Y)^3}=0
\]
We have further from \eqref{eq:e5} that
\[
\lim_{|h|\to 0}\Big\|\scurl_y{\partial\over\partial y_q}\left({N^r(x+h,\cdot)-N^r(x,\cdot)\over h}-\chi^r(x,\cdot)\right)\Big\|_{L^2(Y)^3}=0
\]
so
\beqas
\lim_{|h|\to 0}\Big\|{\partial\over\partial y_p}\left({N^r(x+h,\cdot)-N^r(x,\cdot)\over h}-\chi^r(x,\cdot)\right)\Big\|_{H^1(Y)^3}=0.
\eeqas
Therefore 
\[
\lim_{|h|\to 0}\Big\|{N^r(x+h,\cdot)-N^r(x,\cdot)\over h}-\chi^r(x,\cdot)\Big\|_{H^2(Y)^3}=0. 
\]
This shows that $N^r\in C^1(\bar D,H^2(Y))^3\subset C^1(\bar D,C(\bar Y))^3$. 
\eproof
\begin{proposition}
Assume that $b(x,y)\in C^1(\bar D,C^2(\bar Y))^{3\times 3}$. The solution $w^r$ of cell problem \eqref{eq:cellprobw} belongs to $C^1(\bar D, C^1(\bar Y))$.
\end{proposition}
\bproof
The cell problem \eqref{eq:cellprobw} can be written as
\[
-\nabla_y\cdot(b(x,y)\nabla_yw^r(x,y))=\nabla_y(b(x,y)e_r).
\]
Fixing $x\in \bar D$, the right hand side is bounded uniformly in $H^1(Y)$ so $w^r(x,\cdot)$ is uniformly bounded in $H^3(Y)$ from elliptic regularity (see McLean \cite{McLean} Theorem 4.16).  
For $h\in \IR^3$, we note that
\beqas
-\nabla_y\cdot[b(x,y)\nabla_y(w^r(x+h,y)-w^r(x,y))]=\nabla_y\cdot[(b(x+h,y)-b(x,y))e_r]\\
+\nabla_y\cdot[(b(x+h,y)-b(x,y))\nabla_yw^r(x+h,y)]:=i_1.
\eeqas
As $\int_Yw^r(x,y)dy=0$, we have
\beqas
&&\|w^r(x+h,y)-w^r(x,y)\|_{H^1(Y)}\le c\|\nabla_y(w^r(x+h,y)-w^r(x,y))\|_{L^2(Y)}\\
&&\le c\|(b(x+h,y)-b(x,y))e_r\|_{L^2(Y)}+c\|(b(x+h,y)-b(x,y))\nabla_yw^r(x+h,y)\|_{L^2(Y)}
\eeqas
which converges to 0 when $|h|\to 0$. 
Fixing $x\in\bar D$, we then have from Theorem 4.16 of \cite{McLean} that
\be
\|w^r(x+h,y)-w^r(x,y)\|_{H^3(Y)}\le \|w^r(x+h,y)-w^r(x,y)\|_{H^1(Y)}+\|i_1\|_{H^1(Y)}
\label{eq:w1}
\ee
which converges to 0 when $|h|\to 0$.
Fixing an index $q$, let $h\in \IR^3$ be a vector whose components are all zero except the $q$th component. Let $\eta(x,\cdot)\in H^1_\#(Y)/\IR$ be the solution of the problem
\[
-\nabla_y\cdot[b(x,y)\nabla_y\eta(x,y)]=\nabla_y\cdot\left[{\partial b\over\partial x_q}e_r\right]+\nabla_y\cdot\left[{\partial b\over\partial x_q}\nabla_yw^r(x,y)\right].
\]
We have
\beqas
-\nabla_y\cdot\left[b(x,y)\nabla_y\left({w^r(x+h,y)-w^r(x,y)\over h}-\eta\right)\right]=\nabla_y\cdot\left[\left({b(x+h,y)-b(x,y)\over h}-{\partial b\over\partial x_q}\right)e_r\right]\\
+\nabla_y\cdot\left[\left({b(x+h,y)-b(x,y)\over h}-{\partial b\over\partial x_q}(x,y)\right)\nabla_yw^r(x+h,y)\right]\\
+\nabla_y\left[{\partial b(x,y)\over\partial x_q}\left(\nabla_yw^r(x+h,y)-\nabla_yw^r(x,y)\right)\right]:=i_2.
\eeqas
From \eqref{eq:w1} and the regularity of $b$, $\lim_{|h|\to 0}\|i_2(x,\cdot)\|_{H^1(Y)}=0$. As $\int_Yw^r(x,y)dy=0$ and $\int_Y\eta(x,y)dy=0$, we have that
\[
\lim_{|h|\to 0}\Big\|{w^r(x+h,\cdot)-w^r(x,\cdot)\over h}-\eta(x,\cdot)\Big\|_{H^1(Y)}=0.
\]
Therefore from Theorem 4.16 of \cite{McLean}, we have
\[
\Big\|{w^r(x+h,\cdot)-w^r(x,\cdot)\over h}-\eta(x,\cdot)\Big\|_{H^3(Y)}\le\Big\|{w^r(x+h,\cdot)-w^r(x,\cdot)\over h}-\eta(x,\cdot)\Big\|_{H^1(Y)}+\|i_2\|_{H^1(Y)} 
\]
which converges to 0 when $|h|\to 0$. 
Thus $w^r\in C^1(\bar D,H^3(Y))\subset C^1(\bar D,C^1(\bar Y))$. 
\eproof

For the regularity of the solution $u_0$ of the homogenized problem \eqref{eq:homprob} we have the following result.
\begin{proposition}
Assume that $D$ is a Lipschitz polygonal domain, and the coefficient $a(\cdot,y)$ is uniformly Lipschitz with respect to $x$, then there is a constant $0<s<1$ so that $\curl u_0\in H^s(D)$.
\end{proposition}
\bproof
When $a(x,y)$ is uniformly Lipschitz with respect to $x$, from \eqref{eq:cellprobN}, $\|\scurl_y N^r(x,\cdot)\|_{L^2(Y)}$ is a Lipschitz function of $x$, so from  \eqref{eq:a0} we have that $a^0$ is  Lipschitz with respect to $x$. As $a^0$ is positive definite, $(a^0)^{-1}$ is Lipschitz. Let $U=a^0\curl u_0$. We have from \eqref{eq:homprob} that $U\in H(\curl;D)$, ${\rm div}((a^0)^{-1}U)=0$ and $(a^0)^{-1}U\cdot n=0$ on $\partial D$ where $n$ is the outward normal vector on $\partial D$. The conclusion follows from Lemma 4.2 of Hiptmair \cite{Hiptmair}.
\eproof
\begin{remark}
If $a^0$ is isotropic, we have from \eqref{eq:homprob} that
\[
\curl\curl u_0 =-{1\over a^0}\nabla a^0\times\curl u_0-{1\over a^0}b^0u_0+{1\over a^0}f\in L^2(D)
\]
so $u_0\in H^1(\scurl,D)$. However, even if $a$ is isotropic, $a^0$ may not be isotropic.
\end{remark}
\bibliographystyle{plain}
\bibliography{homerror}

\end{document}